\documentclass[11pt]{amsart}

\usepackage{latexsym,amssymb,amsmath,graphics}
\usepackage[dvips]{graphicx}

\def\pp{{\mathcal P}}

\def\ss{{\mathcal S}}

\def\ffi{\varphi}
\def\eps{\varepsilon}
\def\dst{\displaystyle}

\def\C{{\mathbb{C}}}

\def\K{{\mathbb{K}}}

\def\P{{\mathbb{P}}}

\def\R{{\mathbb{R}}}
\def\S{{\mathbb{S}}}

\def\Z{{\mathbb{Z}}}

\newcommand{\norm}[1]{{\left\|{#1}\right\|}}
\newcommand{\ent}[1]{{\left[{#1}\right]}}
\newcommand{\abs}[1]{{\left|{#1}\right|}}
\newcommand{\scal}[1]{{\left\langle{#1}\right\rangle}}

\newenvironment{definition}[1][]{\vskip3pt\noindent\sl\textbf{Definition.}\ }{\rm\vskip3pt}

\newtheorem{lemma}{Lemma}[section]
\newtheorem{proposition}[lemma]{Proposition}
\newtheorem{theorem}[lemma]{Theorem}
\newtheorem{corollary}[lemma]{Corollary}

\newtheorem{examplenum}[lemma]{Example}

\begin{document}

\title{Uncertainty principles for orthonormal sequences}
\author{Philippe Jaming}
\address{Universit\'e d'Orl\'eans\\
Facult\'e des Sciences\\ 
MAPMO\\ BP 6759\\ F 45067 Orl\'eans Cedex 2\\
France}
\email{Philippe.Jaming@univ-orleans.fr}
\thanks{The first author was partially supported by a {European Commission} grant on
Harmonic Analysis and Related Problems 2002-2006 IHP Network
(Contract Number: HPRN-CT-2001-00273 - HARP),
by the Balaton program EPSF, and by the Erwin Schr\"odinger Insitute.}

\author{Alexander M. Powell}
\address{Department of Mathematics\\
1326 Stevenson Center\\
Vanderbilt University\\
Nashville, TN 37240\\
USA}
\email{alexander.m.powell@vanderbilt.edu}
\thanks{The second author was partially supported by NSF DMS Grant 0504924
and by an Erwin Schr\"odinger Institute Junior Research Fellowship}

\begin{abstract}
The aim of this paper is to provide complementary quantitative 
extensions of two 
results
of H.S. Shapiro on the time-frequency concentration of orthonormal sequences in $L^2 (\R)$.
More precisely, Shapiro proved that if the elements of an orthonormal sequence
and their Fourier transforms are all pointwise bounded by a fixed function in $L^2(\R)$
then the sequence is finite.  In a related result, Shapiro also proved that if the elements of
an orthonormal sequence and their Fourier transforms have uniformly bounded means and dispersions then
the sequence is finite.

This paper gives quantitative bounds on the size of the finite orthonormal sequences in Shapiro's 
uncertainty principles. The bounds are obtained by using prolate sphero\"\i dal wave functions 
and combinatorial estimates on the number of elements in a spherical code.  
Extensions for Riesz bases and different measures of time-frequency concentration are also given.

\end{abstract}



\date{}


\maketitle


\section{Introduction}

The uncertainty principle in harmonic analysis is a class
of theorems which state that a nontrivial function 
and its Fourier transform can not both be too sharply localized.
For background on different appropriate notions of localization
and an overview on the recent renewed interest in mathematical formulations
of the uncertainty principle, see the survey \cite{JP.FS}.
This paper will adopt the broader view 
that the uncertainty principle can be seen not only as a statement
about the time-frequency localization of a single function but also as a
statement on the degradation of localization when one considers successive elements
of an orthonormal basis.  In particular, the results that we consider
show that the elements of an orthonormal basis as well as their Fourier transforms
can not be uniformly concentrated in the time-frequency plane. 

Hardy's Uncertainty Principle \cite{JP.Ha} may be viewed as an early theorem of this type.
To set notation, define the 
\emph{Fourier transform} of $f\in L^1(\R)$ by
$$
\widehat{f}(\xi)=\int f(t)e^{-2i\pi t\xi} {d}t,
$$
and then extend to $L^2(\R)$ in the usual way. 

\medskip

\begin{theorem} [Hardy's Uncertainty Principle]   Let $a,b,C,N>0$ be positive real numbers and let
$f \in L^2(\R)$. Assume that for almost every $x,\xi \in \R$,
\begin{equation} \label{eq_hardy}
|f(x)|\leq C(1+|x|)^Ne^{-\pi a|x|^2}\ \ \  \hbox{ and } \ \ \
|\widehat{f}(\xi)|\leq C(1+|\xi|)^Ne^{-\pi b|\xi|^2}.
\end{equation}
The following hold:
\begin{itemize}
\item If $ab>1$ then $f=0$.
\item If $ab=1$ then $f(x)=P(x)e^{-\pi a|x|^2}$ for some polynomial $P$ of degree at most $N$.
\end{itemize}
\end{theorem}

This theorem has been further generalized where the pointwise condition (\ref{eq_hardy})
is replaced by
integral conditions in \cite{JP.BDJ}, and by distributional conditions in \cite{JP.De}.
Also see \cite{AP.GZ} and \cite{AP.HL}.
One may interpret Hardy's theorem by saying that the set of functions which, along their
Fourier transforms, is bounded by $C(1+|x|)^Ne^{-\pi|x|^2}$ is finite dimensional,
in the sense that its
span is a finite dimensional subspace of $L^2(\R)$.

In the case $ab<1$, the class of functions satisfying the condition (\ref{eq_hardy}) 
has been fully described by B. Demange \cite{JP.De}. 
In particular, it is an infinite dimensional subset of $L^2(\R)$.
Nevertheless, it can not contain an infinite orthonormal sequence.
Indeed, this was first proved by Shapiro in \cite{JP.Sh}:

\begin{theorem} [Shapiro's Umbrella Theorem] \label{ap.sh-umbrel}
Let $\ffi,\psi \in  L^2(\R)$.  If
$\{ e_k \} \subset L^2 (\R)$ is an orthonormal sequence of functions such that for all $k$ and for almost all $x,\xi\in\R$,
$$
|e_k(x)|\leq |\ffi(x)|\quad\mbox{and}\quad|\widehat{e}_k(\xi)|\leq |\psi(\xi)|,
$$
then the sequence $\{e_k\}$ is finite.
\end{theorem}

Recent work of A. De Roton, B. Saffari, H.S. Shapiro, G. Tennenbaum,
see \cite{JP.DeR}, shows that the assumption $\ffi,\psi \in L^2 (\R)$ can not be substantially weakened.
Shapiro's elegant proof of Theorem \ref{ap.sh-umbrel} uses a compactness argument of Kolmogorov, 
see \cite{JP.Sh2},
but does not give a bound on the number of elements in the finite sequence.

A second problem of a similar nature studied by Shapiro in \cite{JP.Sh} is that
of bounding the means and variances of orthonormal sequences.
For $f\in L^2(\R)$ with $\norm{f}_2=1$, we define the following associated
\emph{mean}
$$
\mu(f)={\rm Mean}(|f|^2) = \int t|f(t)|^2\mbox{d}t,
$$
and the associated \emph{variance}
$$
\Delta^2(f)= {\rm Var} (|f|^2) = \int|t-\mu(f)|^2|f(t)|^2\mbox{d}t.
$$
It will be convenient to work also with the \emph{dispersion} $\Delta(f)\equiv\sqrt{\Delta^2(f)}$.
In \cite{JP.Sh}, Shapiro posed the question of determining for which sequences of real numbers
$\{a_n\}_{n=0}^{\infty},$ $\{b_n\}_{n=0}^{\infty},$ $\{c_n\}_{n=0}^{\infty},$
$\{d_n\}_{n=0}^{\infty} \subset \R$ there exists an orthonormal basis
$\{e_n\}_{n=0}^{\infty}$ for $L^2 (\R)$ such that for all $n \geq 0$
$$
\mu(e_n)=a_n,\ \mu(\widehat{e}_n)=b_n,\ \Delta(e_n)=c_n,\ \Delta(\widehat{e}_n)=d_n.
$$
Using again Kolmogorov's compactness argument, 
he proved the following, \cite{JP.Sh}:

\begin{theorem} [Shapiro's Mean-Dispersion Principle]\label{ap.sh-meandisp}
There does not exist an infinite orthonormal sequence $\{e_n\}_{n=0}^{\infty}
\subset L^2(\R)$ such that
all four of $\mu(e_n)$, $\mu(\widehat{e}_n)$, $\Delta(e_n), \Delta(\widehat{e}_n)$ are uniformly bounded.
\end{theorem}

\medskip

An extension of this theorem in \cite{JP.Po} shows that if
$\{e_n\}_{n=0}^{\infty}$ is an orthonormal basis for $L^2(\R)$ then
two dispersions and one mean $\Delta(e_n), \Delta(\widehat{e}_n),\mu(e_n)$
can not all be uniformly bounded.
Shapiro recently pointed out a nice alternate proof of this result using
the Kolmogorov compactness theorem from \cite{JP.Sh}.  The case for two means and one dispersion 
is different.   In fact, it is possible to construct an 
orthonormal basis $\{e_n\}_{n=0}^{\infty}$ for $L^2(\R)$ such that
the two means and one dispersion $\mu(e_n)$, $\mu(\widehat{e}_n),\Delta(e_n)$ 
are uniformly bounded, see \cite{JP.Po}.  

Although our focus will be on Shapiro's theorems,
let us also briefly refer the reader to some other work in the literature concerning
uncertainty principles for bases.
The classical Balian-Low theorem states that if a set of lattice coherent states forms
an orthonormal basis for $L^2(\R)$ then the window function satisfies a strong
version of the uncertainty principle, e.g., see \cite{AP.CP,AP.GHHK}.  
For an analogue concerning dyadic orthonormal wavelets, see \cite{AP.battle}.

\subsection*{Overview and main results}
The goal of this paper is to provide quantitative versions of Shapiro's
Mean-Dispersion Principle and Umbrella
Theorem, i.e., Theorems \ref{ap.sh-umbrel} and \ref{ap.sh-meandisp}.

Section \ref{sec2.ap} addresses the Mean-Dispersion Theorem.
The main results of this section are contained in 
Section \ref{ap.sharp.md.sec} where we prove a sharp quantitative 
version of Shapiro's Mean-Dispersion Principle.  
This result is sharp, but the method of proof
is not easily applicable to more general
versions of the problem.  
Sections \ref{ap.herm.sec} and \ref{ap.rayleigh.sec} respectively 
contain necessary background on Hermite functions and the Rayleigh-Ritz technique which is needed
in the proofs.  Section \ref{ap.riesz.shmd.sec} proves a version of the mean-dispersion theorem for Riesz bases. 

Section \ref{sec3.ap} addresses the Umbrella Theorem
and variants of the Mean-Dispersion Theorem.
The main results of this section are contained in 
Section \ref{gen.md.ap} where we prove a 
quantitative version of the Mean-Dispersion Principle for a generalized
notion of dispersion, and in Section \ref{JP.sec:QUT} where
we prove a quantitative version of Shapiro's Umbrella Theorem. 
Explicit bounds on the size of possible orthonormal sequences are given in particular cases.  
Since the methods of Section \ref{sec2.ap} are no longer
easily applicable here, we adopt an approach based on geometric combinatorics.
Our results use estimates on the
size of spherical codes, and the theory of prolate sphero\"\i dal wavefunctions.
Section \ref{ap.sphcode.sec} contains background results on spherical codes, 
including the Delsarte, Goethals, Seidel bound.
Section \ref{ap.apponb.sec} proves some necessary results on projections of one set of orthonormal functions onto
another set of orthonormal functions.  Section \ref{ap.psw.sec} gives an overview of the prolate sphero\"\i dal wavefunctions
and makes a connection between projections of orthonormal functions and spherical codes.
Section \ref{ap.rieszang.sec} concludes with extensions to Riesz bases.

\section{Growth of means and dispersions} \label{sec2.ap}

In this section, we use the classical Rayleigh-Ritz technique to give a quantitative version of
Shapiro's Mean-Dispersion Theorem. We also prove that, in this sense, the Hermite
basis is the best concentrated orthonormal basis of $L^2(\R)$.

\subsection{The Hermite basis} \label{ap.herm.sec}

Results of this section can be found in \cite{JP.FS}.
The Hermite functions are defined by
$$
h_k(t)=\frac{2^{1/4}}{\sqrt{k!}}\left(-\frac{1}{\sqrt{2\pi}}\right)^ke^{\pi t^2} \left(\frac{{d}}{{d}t}\right)^ke^{-2\pi t^2}.
$$
It is well known that the Hermite functions are eigenfunctions of the Fourier transform,
satisfy $\widehat{h_k}=i^{-k}h_k$, and 
form an orthonormal basis for $L^2(\R)$.
Let us define the {\em Hermite operator} $H$ for
functions $f$ in the Schwartz class by
$$Hf(t) = -\frac{1}{4\pi^2}\frac{d^2}{dt^2}f(t) + t^2 f(t).$$
It is easy to show that 
\begin{equation}
\label{JP.eigen_val}
Hh_k = \left( \frac{2k+1}{2\pi} \right) h_k,
\end{equation}
so that $H$ may also be seen as the densely defined,
positive, self-adjoint, unbounded operator on $L^2(\mathbb{R})$ defined by
$$
Hf = \sum_{k=0}^{\infty} \frac{2k+1}{2\pi} \langle f, h_k \rangle h_k.
$$
From this, it immediately follows that, for each $f$ in the domain of $H$ 
\begin{eqnarray}
\langle Hf,f \rangle&=& 
\sum_{k=0}^{\infty} \frac{2k+1}{2\pi} |\langle f, h_k \rangle |^2
= \int |t|^2 |f(t)|^2 dt + \int |\xi|^2 |\widehat{f}(\xi)|^2 d\xi \label{JP.eq:fundherm} \\
&=&\mu(f)^2\norm{f}_2^2+\Delta^2(f)+\mu(\widehat{f})^2\|\widehat{f}\|_2^2+\Delta^2(\widehat{f}).
\notag 
\end{eqnarray}

\subsection{The Rayleigh-Ritz Technique} \label{ap.rayleigh.sec}

The Rayleigh-Ritz technique
is a useful tool for estimating eigenvalues of operators, see \cite[Theorem XIII.3, page 82]{JP.RS4}.

\begin{theorem} [The Rayleigh-Ritz Technique]
\label{JP.ral_ritz_thm}
Let $H$ be a positive self-adjoint operator and define
$$
\lambda_k (H) = \sup_{\varphi_0, \cdots, \varphi_{k-1}} \ \ 
\inf_{\psi \in [\varphi_0, \cdots, \varphi_{k-1}]^{\perp},\|\psi\|_2 \leq 1, 
\psi \in D(H)} 
\langle H \psi, \psi \rangle.$$
where $D(H)$ is the domain of $H$. Let $V$ be a $n+1$ dimensional
subspace, $V\subset D(H),$ and let $P$ be the orthogonal projection onto $V$.
Let $H_V = PHP$ and let $\widetilde{H_V}$ denote the restriction of
$H_V$ to $V$.  
Let $\mu_0 \leq \mu_1 \leq \cdots \leq \mu_n$ be the eigenvalues of 
$\widetilde{H_V}$
Then
$$
\lambda_k (H) \leq \mu_k, \ \ \ k=0, \cdots, n.
$$
\end{theorem}

The following corollary is a standard and useful application
of the Rayleigh-Ritz technique.  For example, \cite[Chapter 12]{JP.LL}
contains a version in the setting of Schr\"odinger operators.

\begin{corollary}
\label{JP.trace_cor}
Let $H$ be a positive self-adjoint operator, and
let $\varphi_0, \cdots, \varphi_n$ be an orthonormal set of functions.  Then
\begin{equation}
\label{JP.eq:trace}
\sum_{k=0}^n \lambda_k(H)  \leq \sum_{k=0}^n
\langle H\varphi_k, \varphi_k \rangle.
\end{equation}
\end{corollary}

\begin{proof}
If some $\varphi_k \notin D(H)$ then positivity of $H$ implies that \eqref{JP.eq:trace} trivially
holds since the right hand side of the equation would be infinite.
We may thus assume that $\varphi_0, \cdots, \varphi_n\in D(H)$.

Define the $n+1$ dimensional subspace $V = {\rm span} \ \{\varphi_k\}_{k=0}^n$
and note that the operator $\widetilde{H_V}$
is given by the matrix
$M =[\langle H \varphi_j, \varphi_k \rangle]_{0\leq j,k \leq n}$.
Let $\mu_0, \cdots, \mu_n$ be the eigenvalues of $\widetilde{H_V}$,
i.e., of the matrix $M$.
By Theorem \ref{JP.ral_ritz_thm},
$$
\sum_{k=0}^n \lambda_k (H) \leq \sum_{k=0}^n \mu_k = {\rm Trace} (M)
= \sum_{k=0}^n \langle H\varphi_k, \varphi_k \rangle
$$
which completes the proof of the corollary.
\end{proof}

\subsection{The Sharp Mean-Dispersion Principle} \label{ap.sharp.md.sec}

\begin{theorem} [Mean-Dispersion Principle]
\label{JP.th:optherm}
Let $\{e_k\}_{k=0}^{\infty}$ be any orthonormal sequence in $L^2(\R)$.
Then for all $n\geq 0$,
\begin{equation}
\label{JP.eq:optherm}
\sum_{k=0}^n \left( \Delta^2(e_k) + \Delta^2(\widehat{e_k}) + 
|\mu(e_k)|^2 + |\mu(\widehat{e_k})|^2 \right) \geq 
\frac{(n+1)(2n+1)}{4\pi}.
\end{equation}
Moreover, if equality holds for all $0 \leq n\leq n_0$, then there exists $\{c_k\}_{n=0}^{n_0}
\subset \C$ such that $|c_k|=1$ and $e_k=c_kh_k$ for each $0 \leq k \leq n_0$.
\end{theorem}

\begin{proof}
Since $H$ is positive and self-adjoint, one may use Corollary 
\ref{JP.trace_cor}.  It follows from Corollary \ref{JP.trace_cor} that for each
$n \geq 0$ one has
\begin{equation}
\label{JP.trace_eq}
\sum_{k=0}^n \frac{2k+1}{2\pi} \leq 
\sum_{k=0}^n \langle He_k, e_k \rangle.
\end{equation}
From (\ref{JP.eq:fundherm}), note that since $\|e_k\|_2=1$,
$$
\langle H e_k, e_k \rangle  = 
\Delta^2 (e_k) + \Delta^2 (\widehat{e_k}) + |\mu(e_k)|^2 
+ |\mu(\widehat{e_k})|^2.
$$
This completes the proof of the first part.

Assume equality holds in (\ref{JP.eq:optherm}) for all $n=0,\ldots,n_0$,
in other terms that, for $n=0,\ldots,n_0$,
$$
\scal{He_n,e_n}=
\Delta^2(e_n) + \Delta^2(\widehat{e_n}) + 
|\mu(e_n)|^2 + |\mu(\widehat{e_n})|^2 =
\frac{2n+1}{2\pi}.
$$
Let us first apply (\ref{JP.eq:fundherm}) for $f=e_0$:
$$
\sum_{k=0}^{\infty} \frac{2k+1}{2\pi} |\langle e_0, h_k \rangle |^2
=\scal{He_0,e_0}=\frac{1}{2\pi}
=\sum_{k=0}^{\infty} \frac{1}{2\pi} |\langle e_0, h_k \rangle |^2
$$
since $\|e_0\|_2=1$.
Thus, for $k\geq1$, one has $\scal{e_0,h_k}=0$ and hence $e_0=c_0h_0$. 
Also $\|e_0\|_2=1$ implies $|c_0|=1$.
Next, assume that we have proved $e_k=c_kh_k$
for $k=0,\ldots,n-1$.
Since $e_n$ is orthogonal to $e_k$ for $k<n,$ one has $\scal{e_n,h_k}=0$. 
Applying (\ref{JP.eq:fundherm}) for $f=e_n$
we obtain that,
$$
\sum_{k=n}^{\infty} \frac{2k+1}{2\pi} |\langle e_n, h_k \rangle |^2
=\sum_{k=0}^{\infty} \frac{2k+1}{2\pi} |\langle e_n, h_k \rangle |^2
=\scal{He_n,e_n}=\frac{2n+1}{2\pi}
=\sum_{k=n}^{\infty} \frac{2n+1}{2\pi} |\langle e_n, h_k \rangle |^2.
$$
Thus $\scal{e_n,h_k}=0$ for $k>n$. It follows that $e_n=c_nh_n$.
\end{proof}

\begin{examplenum}
For all $n \geq 0$, the Hermite functions satisfy 
$$\mu(h_n) = \mu(\widehat{h_n}) =0 \ \ \ \hbox{ and } \ \ \
\Delta^2 (h_n) =\Delta^2 (\widehat{h_n}) = \frac{2n+1}{4\pi}.$$
\end{examplenum}
For comparison, let us remark that 
Bourgain has constructed an orthonormal basis $\{b_n\}_{n=1}^{\infty}$
for $L^2(\R)$, see \cite{JP.B}, 
which 
satisfies $\Delta^2 (b_n) \leq \frac{1}{2\sqrt{\pi}} + \eps$ 
and $\Delta^2 (\widehat{b_n}) \leq \frac{1}{2\sqrt{\pi}} + \eps$.  However,
it is difficult to control the growth of $\mu(b_n), \mu(\widehat{b_n})$ in this
construction.
For other bases with more structure, see the
related work in \cite{AP.BCGP} that constructs 
an orthonormal basis of lattice coherent states $\{g_{m,n}\}_{m,n\in\Z}$
for $L^2(\R)$  
which is logarithmically close to having 
uniformly bounded dispersions.  The means 
$(\mu(g_{m,n}), \mu( \widehat{g_{m,n}}))$ for this basis lie on a translate of
the lattice $\Z \times \Z$.


It is interesting to note that if one takes $n=0$ in Theorem \ref{JP.th:optherm} then this yields the usual 
form of Heisenberg's uncertainty principle (see \cite{JP.FS} for
equivalences between uncertainty principles with sums and products).
In fact, using (\ref{JP.eq:fundherm}), Theorem \ref{JP.th:optherm}
also implies a more general version of Heisenberg's uncertainty
principle that is implicit in \cite{JP.FS}.
In particular, if $f\in L^2(\R)$ with $\norm{f}_2=1$ is orthogonal to $h_0,\ldots,h_{n-1}$ then
$$
\Delta^2(f) + \Delta^2(\widehat{f}) + 
|\mu(f)|^2 + |\mu(\widehat{f})|^2 \geq \frac{2n+1}{2\pi}.
$$

For instance, if $f$ is odd, then $f$ is orthogonal to $h_0$,
and $\mu(f)=\mu(\widehat{f})=0$. Using the usual scaling trick,
we thus get the well known fact that
the optimal constant in Heisenberg's inequality, e.g., see \cite{JP.FS}, is given as follows
$$
\Delta(f)\Delta(\widehat{f}) \geq
\begin{cases}
\frac{1}{4\pi}\norm{f}_2^2&\mbox{in general}\\
\frac{3}{4\pi}\norm{f}_2^2&\mbox{if $f$ is odd}
\end{cases}.
$$

\begin{corollary} \label{ap.cor.opt.md}
Fix $A>0.$  If $\{e_k\}_{k=0}^{n} \subset L^2(\R)$ is an orthonormal sequence and
for $k=0,\ldots,n$, satisfies
$$
|\mu(e_k)|,\ |\mu(\widehat{e_k})|,\ 
\Delta(e_k),\ \Delta(\widehat{e_k})\leq A,
$$ 
then
$n\leq 8\pi A^2$.
\end{corollary}

\begin{proof}
According to Theorem \ref{JP.th:optherm}
$$
4(n+1)A^2\geq
\sum_{k=0}^n \left( \Delta^2(e_k) + \Delta^2(\widehat{e_k}) + 
|\mu(e_k)|^2 + |\mu(\widehat{e_k})|^2 \right) \geq \frac{(n+1)(2n+1)}{4\pi}.
$$
It follows that $2n+1\leq 16\pi A^2$.
\end{proof}

This may also be stated as follows:

\begin{corollary}
\label{JP.cor:cor2}
If $\{e_k\}_{k=0}^{\infty} \subset L^2(\R)$ is an orthonormal sequence, then
for every $n$,
$$
\max\{|\mu(e_k)|,\ |\mu(\widehat{e_k})|,\ \Delta(e_k),\ \Delta(\widehat{e_k})\ : 0 \leq k \leq n\}
\geq \sqrt{\frac{2n+1}{16\pi}}.
$$
\end{corollary}

\subsection{An extension to Riesz bases} \label{ap.riesz.shmd.sec}

Recall that $\{x_k\}_{k=0}^{\infty}$ is a
\emph{Riesz basis} for $L^2(\R)$ if there 
exists an isomorphism, $U:L^2(\R)\to L^2(\R)$, called the \emph{orthogonalizer} 
of $\{x_k\}_{k=0}^{\infty}$, such that
$\{Ux_k\}_{k=0}^{\infty}$ is an orthonormal basis for $L^2(\R)$. 
It then follows that, for every 
$\{a_n\}_{n=0}^{\infty} \in\ell^2$,
\begin{equation}
\frac{1}{\norm{U}^2}\sum_{n=0}^{\infty} |a_n|^2
\leq \norm{\sum_{n=0}^{\infty} a_nx_n}_{2}^2\leq\norm{U^{-1}}^2\sum_{n= 0}^{\infty} |a_n|^2.
\label{JP.eq:riesz}
\end{equation}

One can adapt the results of the previous sections to Riesz bases.
To start, note that the Rayleigh-Ritz technique leads to the 
following, cf. \cite[Theorem XIII.3, page 82]{JP.RS4}:

\begin{lemma}
\label{JP.lem:traceR}
Let $H$ be a positive, self-adjoint, densely defined operator on $L^2 (\R)$, 
and let $\{x_k\}_{k=0}^{\infty}$ be a Riesz basis for $L^2(\R)$ with 
orthonormalizer $U$. Then, for every $n\geq 0$,
\begin{equation}
\label{JP.eq:traceR}
\sum_{k=0}^n \lambda_k(H)  \leq \norm{U}^2 \sum_{k=0}^n
\langle Hx_k, x_k \rangle.
\end{equation}
\end{lemma}

\begin{proof} Let us take the notations of the proof of Corollary \ref{JP.trace_cor}.
Write $\ffi_k=Ux_k$, it is then enough to notice that
$$
M=[\scal{Hx_k,x_k}]=[\scal{HU^{-1}x_k,U^{-1}x_k}]=[\scal{U^{-1}\,^*HU^{-1}x_k,x_k}].
$$
As $U^{-1}\,^*HU^{-1}$ is a positive operator, that the Rayleigh-Ritz theorem gives
$$
\sum_{k=0}^n\langle Hx_k, x_k \rangle\geq \sum_{k=0}^n\lambda_k(U^{-1}\,^*HU^{-1}).
$$
But,
\begin{eqnarray*}
\lambda_k(U^{-1}\,^*HU^{-1})&=&\sup_{\varphi_0, \cdots, \varphi_{k-1}} \ \ 
\inf_{\psi \in [\varphi_0, \cdots, \varphi_{k-1}]^{\perp},\|\psi\|_2\leq1}
\scal{U^{-1}\,^*HU^{-1}\psi,\psi}\\
&=& \sup_{\varphi_0, \cdots, \varphi_{k-1}}
\inf_{\psi \in [\varphi_0, \cdots, \varphi_{k-1}]^{\perp},\|\psi\|_2\leq1}
\scal{HU^{-1}\psi,U^{-1}\psi}\\
&=& \sup_{\varphi_0, \cdots, \varphi_{k-1}}
\inf_{\tilde\psi \in [U^*\varphi_0, \cdots, U^*\varphi_{k-1}]^{\perp},\|U\tilde\psi\|_2\leq1}
\scal{H\tilde\psi,\tilde\psi}
\end{eqnarray*}
and, as $\|U\tilde\psi\|_2\leq \norm{U}\,\|\tilde\psi\|_2$,
$$
\lambda_k(U^{-1}\,^*HU^{-1})\geq \frac{1}{\norm{U}^2}
\sup_{\tilde\varphi_0, \cdots, \tilde\varphi_{k-1}}
\inf_{\tilde\psi \in [\tilde\varphi_0, \cdots, \tilde\varphi_{k-1}]^{\perp},\|\tilde\psi\|_2\leq1}
\scal{H\tilde\psi,\tilde\psi}=\frac{1}{\norm{U}^2}\lambda_k(H).
$$
\end{proof}

Adapting the proofs of the previous section, we obtain the following corollary.

\begin{corollary}
If $\{ x_k \}_{k=0}^{\infty}$ is a Riesz basis for $L^2(\R)$ with orthonormalizer $U$ then
for all $n$,
$$
\sum_{k=0}^n \left( \Delta^2(x_k) + \Delta^2(\widehat{x_k}) + 
|\mu(x_k)|^2 + |\mu(\widehat{x_k})|^2 \right) \geq 
\frac{(n+1)(2n+1)}{4\pi\norm{U}^2}.
$$
Thus, for every $A>0$, there are at most $8\pi A^2\norm{U}^2$ elements of the
basis $\{ x_k \}_{k=0}^{\infty}$ 
such that $|\mu(e_n)|$, $|\mu(\widehat{e_n})|$, $\Delta(e_n)$, $\Delta(\widehat{e_n})$
are all bounded by $A$. In particular,
$$
\max\{|\mu(x_k)|,\ |\mu(\widehat{x_k})|,\ \Delta(x_k),\ \Delta(\widehat{x_k}):\ 0 \leq k \leq n \}
\geq \frac{1}{\norm{U}}\sqrt{\frac{2n+1}{16\pi}}.
$$
\end{corollary}

\section{Finite dimensional approximations, spherical codes and the Umbrella Theorem} \label{sec3.ap}

\subsection{Spherical codes} \label{ap.sphcode.sec}

Let $\K$ be either $\R$ or $\C$, and let $d\geq1$ be a fixed integer. We equip $\K^d$ with the standard
Euclidean scalar product and norm. We denote by $\S_d$ the unit sphere of $\K^d$.

\begin{definition}
Let $A$ be a subset of $\{z\in\K\,:\ |z|\leq 1\}$. A spherical $A$-code is a finite subset $V \subset \S_d$
such that if $u,v \in V$ and $u\not=v$ then $\scal{u,v}\in A$.  

\end{definition}

Let $N_\K(A,d)$ denote the maximal cardinality of a spherical $A$-code.
This notion has been introduced in \cite{JP.DeGoSi} in the case $\K=\R$ where upper-bounds
on $N_\R(A,d)$ have been obtained.  These are important quantities in geometric combinatorics,
and there is a large associated literature.  Apart from \cite{JP.DeGoSi},
the results we use can all be found in \cite{JP.CoSl}.

Our prime interest is in the quantity
$$N_\K^s(\alpha,d)=\begin{cases}N_\R([-\alpha,\alpha],d),&\mbox{when }\K=\R\\
N_\C(\{z\in \C:|z|\leq \alpha\},d),&\mbox{when }\K=\C\\ \end{cases}$$ for $\alpha\in (0,1]$.
Of course $N_\R^s(\alpha,d)\leq N_\C^s(\alpha,d)$.
Using the standard identification of $\C^d$ with $\R^{2d}$, namely
identifying $Z=(x_1+iy_1,\ldots, x_d + iy_d)\in\C^d$ with 
$\tilde Z=(x_1,y_1,\ldots, x_d, y_d)\in\R^{2d}$, we have
$\scal{\tilde Z,\tilde Z'}_{\R^{2d}}=\mbox{Re}\,\scal {Z,Z'}_{\C^d}$.
Thus $N_\C^s(\alpha,d)\leq N_\R^s(\alpha,2d)$.

In dimensions $d=1$ and $d=2$ one can compute the following values for $N_\K^s(\alpha,d)$:
\begin{itemize}
\item $N_\R^s(\alpha,1)=1$
\item If $0 \leq \alpha < 1/2$ then $N_\R^s(\alpha,2)=2$
\item If $\cos\frac{\pi}{N}\leq \alpha<\cos\frac{\pi}{N+1}$ and $3 \leq N$
then $N_\R^s(\alpha,2)=N$.
\end{itemize}
In higher dimensions, one has the following result.
\begin{lemma} \label{ap.sphcodebnd.lem}
If $0\leq \alpha<\frac{1}{d}$ then $N_\K^s(\alpha,d)=d$.
\end{lemma} 
\begin{proof}
An orthonormal basis of $\K^d$
is a spherical $[-\alpha,\alpha]$-code so that $N_\K^s(\alpha,d)\geq d$.
For the converse, let $\alpha<1/d$ and assume towards a contradiction that
$w_0,\ldots,w_d$ is a spherical $[-\alpha,\alpha]$-code.
Indeed, let us show that $w_0,\ldots,w_d$ would be linearly independent in $\K^d$.
Suppose that $\dst\sum_{j=0}^{d}\lambda_jw_j=0,$ and without loss of generality
that $|\lambda_j|\leq|\lambda_0|$ for $j=1,\ldots,d$. 
Then $\lambda_0\norm{w_0}^2=\dst-\sum_{j=1}^{d}\lambda_j\scal{w_j,w_0}$
so that
$|\lambda_0|\leq |\lambda_0|d\alpha.$
As $d\alpha<1$ we get that $\lambda_0=0$ and then $\lambda_j=0$  for all $j$.
\end{proof}

In general, it is difficult to compute $N_\K^s(\alpha,k)$.  A coarse
estimate using volume counting proceeds as follows.
\begin{lemma} \label{ap.expbnd.lem}
If $0 \leq \alpha <1$ is fixed, then there exist constants
$0<a_1<a_2$ and $0<C$ such that for all $d$ 
$$
\frac{1}{C}e^{a_1 d}\leq N_\K^s(\alpha,d)\leq Ce^{a_2 d}.
$$
Moreover, for $\alpha\leq 1/2$ one has 
$N_\K^s(\alpha,d)  \leq 3^d$ if $\K=\R$, and $N_\K^s(\alpha,d) \leq 9^d$ 
if $\K=\C$.
\end{lemma}


\begin{proof}
\noindent The counting argument for the upper bound proceeds as follows.
Let $\{w_j\}_{n=1}^N$ be a spherical $A$-code, with $A =[-\alpha, \alpha]$
or $A= \{z\in \C : |z| \leq \alpha \}.$
For $j\not=k$, one has
$$
\norm{w_j-w_k}^2=\norm{w_j}^2+\norm{w_k}^2+2\mbox{Re}\,\scal{w_j,w_k}
\geq 2-2\alpha.
$$
So, the open balls $\dst B\left(w_j,\sqrt{\frac{1-\alpha}{2}}\right)$ of center $w_j$ and radius
$\dst\sqrt{\frac{1-\alpha}{2}}$ are all disjoint and included in the ball of center $0$ and radius 
$\dst1+\sqrt{\frac{1-\alpha}{2}}$. Therefore
$$
Nc_d\left(\frac{1-\alpha}{2}\right)^{hd/2}\leq c_d\left(1+\sqrt{\frac{1-\alpha}{2}}\right)^{hd}
$$
where $c_d$ is the volume of the unit ball in $\K^d$,
$h=1$ if $\K=\R$ and $h=2$ if $\K=\C$.
This gives the bound $N\leq \left(1+\sqrt{\frac{2}{1-\alpha}}\right)^{hd}$.
Note that for $\alpha\leq 1/2$ we get $N\leq 3^d$ if $\K=\R$ and $N\leq 9^d$ if $\K=\C$.
The lower bound too may be obtained by a volume counting argument, see \cite{JP.CoSl}.
\end{proof}

The work of Delsarte, Goethals, Seidel, \cite{JP.DeGoSi} provides a method for
obtaining more refined estimates on the size of spherical codes.  
For example, taking $\beta = - \alpha$ in Example 4.5 of \cite{JP.DeGoSi}
shows that 
if $\alpha<\frac{1}{\sqrt{d}}$ then
\begin{equation} \label{delsarte_bnd.ap}
\dst N_\R^s(\alpha,d)\leq \frac{(1-\alpha^2)d}{1-\alpha^2d}.
\end{equation}
Equality can only occur for spherical $\{-\alpha,\alpha\}$-codes.
Also, note that if $\alpha=\frac{1}{\sqrt{d}}\sqrt{1-\frac{1}{d^k}}$, then 
$\frac{1-\alpha^2}{1-\alpha^2d}d\sim d^{k+1}$.

\subsection{Approximations of orthonormal bases} \label{ap.apponb.sec}

We now make a connection between the cardinality of spherical codes and 
projections of orthonormal bases.

Let $\mathcal{H}$ be a Hilbert space over $\K$ and let $\Psi=\{\psi_k\}_{k=1}^{\infty}$ be an orthonormal
basis for $\mathcal{H}$. For an integer $d\geq 1$, let $\P_d$ be the orthogonal
projection on the span of $\{\psi_1,\ldots,\psi_d\}$.
For $\eps>0$, we say that an element $f\in \mathcal{H}$ 
is $\eps,d$-approximable
if $\norm{f-\P_df}_{\mathcal{H}}<\eps$, and define 
${\ss}_{\eps,d}$ to be the set of 
all of $f\in {\mathcal{H}}$ with 
$||f||_{\mathcal{H}}=1$ that are $\eps,d$-approximable.
We denote by $A_\K(\eps,d)$ the maximal cardinality of an orthonormal sequence
in ${\ss}_{\eps,d}$.

\begin{examplenum}
Let $\{e_j\}_{j=1}^n$ be the canonical basis for $\R^n$, and let 
$\{\psi_j\}_{j=1}^{n-1}$ be
an orthonormal basis for $V^{\perp}$,
where $V = {\rm span} \{(1,1,\ldots,1)\}$. 
Then $\norm{e_k-P_{n-1}e_k}_2=\frac{1}{\sqrt{n}}$ holds for each
$1 \leq k \leq n$,
and hence $A_\R(\frac{1}{\sqrt{n}},n-1)\geq n$.
\end{examplenum}

Our interest in spherical codes stems from the following result, cf. \cite[Corollary 1]{JP.Po}.

\begin{proposition}
\label{JP.prop:estimate}
If $0<\eps<1/\sqrt{2}$ and  $\alpha=\frac{\eps^2}{1-\eps^2},$
then $A_\K(\eps,d)\leq N_\K^s(\alpha,d)$.
\end{proposition}


\begin{proof}
Let $\{\psi_k\}_{k=1}^{\infty}$ be an orthonormal basis for $\mathcal{H}$,
and let $\ss_{\eps,d}$ and $\P_d$ be as above.
Let $\{f_j\}_{j=1}^N \subset {\mathcal{H}}$ 
be an orthonormal set contained in $\ss_{\eps,d}$.
For each $k=1,\ldots,N$, $j=1,\ldots,d$,
let $a_{k,j}=\scal{f_k,\psi_j}$ and write $\P_df_k:=\dst\sum_{j=1}^da_{k,j}\psi_j$
so that $\norm{f_k-\P_df_k}_{\mathcal{H}}<\eps$.

Write $v_k=(a_{k,1},\ldots,a_{k,d})\in\K^d$ then, for $k\not=l$
\begin{eqnarray}
\scal{v_k,v_l}&=&\scal{\P_df_k,\P_df_l}=\scal{\P_df_k-f_k+f_k,\P_df_l-f_l+f_l}\nonumber\\
&=&\scal{\P_df_k-f_k,\P_df_l-f_l}+\scal{\P_df_k-f_k,f_l}+\scal{f_k,\P_df_l-f_l}\nonumber\\
&=&\scal{\P_df_k-f_k,\P_df_l-f_l}+\scal{\P_df_k-f_k,f_l-\P_df_l}+\scal{f_k-\P_df_k,\P_df_l-f_l}\nonumber\\
&=&\scal{f_k-\P_df_k,\P_df_l-f_l}
\label{JP.eq:orsp}
\end{eqnarray}
since $\P_df_k-f_k$ is orthogonal to $\P_df_l$.
It follows from the Cauchy-Schwarz inequality
that $\abs{\scal{v_k,v_l}}\leq \eps^2$.

On the other hand,
$$
\norm{v_k}_{\K^d}=\norm{\P_df_k}_{\mathcal{H}}
=(\norm{f_k}_{\mathcal{H}}^2-
\norm{f_k-\P_df_k}_{\mathcal{H}}^2)^{1/2}\geq(1-\eps^2)^{1/2}.
$$
It follows that $w_k=\dst\frac{v_k}{\norm{v_k}_{\K^d}}$ 
satisfies, for $k\not=l$,
$\abs{\scal{w_k,w_l}}=
\frac{\abs{\scal{v_k,v_l}}}{\norm{v_k}_{\K^d}\norm{v_l}_{\K^d}}\leq\frac{\eps^2}{1-\eps^2},$
and $\{w_k\}_{k=1}^N$ is a spherical $[-\alpha,\alpha]$-code in $\K^d$.
\end{proof}

Note that the proof only uses orthogonality in a mild way.
Namely, if instead 
$\{f_j\}_{j=1}^N \subset \mathcal{H}$ with $\norm{f_j}_{\mathcal{H}} =1$ 
satisfies
$|\scal{f_j,f_k}| \leq \eta^2$ for $j\neq k$, then
Equation (\ref{JP.eq:orsp}) becomes 
$\scal{v_k,v_l}=\scal{f_k-\P_df_k,\P_df_l-f_l}
+\scal{f_k,f_l},$ so that 
$|\scal{v_k,v_l}| \leq \eps^2+\eta^2$,
and the end of the proof shows that 
$\dst N\leq N_\K^s\left(\frac{\eps^2+\eta^2}{1-\eps^2},d\right)$.

In view of Proposition \ref{JP.prop:estimate}, 
it is natural ask the following question.
Given $\alpha=\frac{\eps^2}{1-\eps^2}$, is there a converse inequality
of the form $N_\K^s(\alpha,d)\leq CA_\K(\eps',d')$
with $C>0$ an absolute constant and $\eps\leq \eps'\leq C\eps$, $d\leq d'\leq Cd$?  Note that for $\eps$ such that $\alpha<1/d$, we have 
$A_\K(\eps,d)=N_\K^s(\alpha,d)=d$.

\subsection{Prolate sphero\"\i dal wave functions} \label{ap.psw.sec}
In order to obtain quantitative versions of Sha\-pi\-ro's theorems, we will
make use of the prolate sphero\"\i dal wave functions.  For a detailed presentation on
prolate sphero\"\i dal wave functions see \cite{JP.SP,JP.LP1,JP.LP2}.

Fix $T,\Omega>0$ and let $\{\psi_n\}_{n=0}^{\infty}$ be the associated prolate spheroidal wave
functions, as defined in \cite{JP.SP}.  $\{\psi_n\}_{n=0}^{\infty}$ is an orthonormal basis for 
$PW_\Omega\equiv 
\{f\in L^2(\R)\,:\ \mbox{supp}\,\widehat{f}\subseteq [-\Omega,\Omega]\},$ 
and the $\psi_n$ 
are eigenfunctions of the differential operator
$$L=(T^2-x^2) \frac{{d^2}}{{d} x^2}
-2x \frac{{d}}{{d}x} - \frac{\Omega^2}{T^2}x^2.$$ 
As in the previous section, for an integer $d\geq 0$, define $\P_d$ to be the 
projection onto
the span of $\psi_0,\ldots,\psi_{d-1},$ and for $\eps>0$ define 
$$
\mathcal{S}_{\eps,d}=\{f\in L^2(\R)\,: \norm{f}_2=1, \norm{f-\P_df}_2 <\eps\}.
$$
For the remainder of the paper, the
orthonormal basis used in the definitions of
$\mathcal{S}_{\eps,d}$, $\P_d$, and 
${A}_{\mathbb{K}} (\eps, d)$, will always be chosen as the prolate
sphero\"\i dal wavefunctions.  Note that these quantities depend on 
the choice of $T,\Omega$.

Finally, let
$$
\pp_{T,\Omega,\eps}=\left\{f\in L^2(\R)\,:\ \int_{|t|>T}|f(t)|^2\,\mbox{d}t\leq\eps^2
\quad\mbox{and}\quad\int_{|\xi|>\Omega}|\widehat{f}(\xi)|^2\,\mbox{d}\xi\leq\eps^2\right\}
$$
and $\pp_{T,\eps}=\pp_{T,T,\eps}$.

\begin{theorem}[Landau-Pollak \cite{JP.LP2}]
\label{JP.th:LPth12}
Let $T,\ \eps$ be positive constants and let $d=\lfloor 4T\Omega\rfloor+1$. Then,
for every $f\in \pp_{T,\Omega,\eps}$,
$$
\norm{f-\P_df}_2^2\leq 49\eps^2\norm{f}_2^2.
$$
In other words, $\pp_{T,\Omega,\eps}\cap\{f\in L^2 (\R)\,:\ \norm{f}_2=1\}\subset\mathcal{S}_{7\eps,d}$.
\end{theorem}

It follows that the first $d=4T^2+1$  elements 
of the prolate sphero\"\i dal basis well approximate
$\pp_{T,\eps}$, and that $\pp_{T,\eps}$ is ``essentially'' $d$-dimensional.

\subsection{Generalized means and dispersions} \label{gen.md.ap}
As an application of the results on prolate spheroidal wavefunctions
and spherical codes, we shall address a more general version of the 
mean-dispersion theorem.  

Consider the following generalized means and variances.
For $p>1$ and $f\in L^2(\R)$ with $\norm{f}_2=1$, we define the following associated
\emph{$p$-variance}
$$
\Delta_p^2(f)=\inf_{a\in\R}\int|t-a|^p|f(t)|^2\mbox{d}t.
$$
One can show that the infimum is actually a minimum and is attained 
for a unique $a\in\R$ that we call the \emph{$p$-mean}
$$
\mu_p(f)= \mbox{arg min}_{a \in \R}\,\int|t-a|^p|f(t)|^2\mbox{d}t.
$$
As before, define the \emph{$p$-dispersion} $\Delta_p(f)\equiv\sqrt{\Delta_p^2(f)}$.

The proof of the Mean-Dispersion Theorem for $p=2$ via the Rayleigh-Ritz technique
relied on the special relation \eqref{JP.eq:fundherm}
of means and dispersions with the Hermite operator.   In general, beyond the case $p=2$,
such simple relations are not present and the techniques of
Section \ref{sec2.ap} are not so easily applicable.  
However, we shall show how to use
the combinatorial techniques from the beginning of this section
to obtain a quantitative version of Theorem \ref{ap.sh-meandisp} for generalized means and dispersions.



The following lemma is a modification of \cite[Lemma 1]{JP.Po}.
\begin{lemma}
\label{JP.lem:Po-prolate}
Let $A>0$ and $p>1$. Suppose $g\in L^2(\R)$, $\norm{g}_2=1$ satisfies
$$
|\mu_p(g)|,\ |\mu_p(\widehat{g})|,\ \Delta_p(g),\ \Delta_p(\widehat{g})\leq A.
$$
Fix $\eps>0$, then $g\in\pp_{A+(A/\eps)^{2/p},\eps}$.
\end{lemma}

This gives a simple proof of a strengthened version of Shapiro's Mean-Dispersion Theorem:
\begin{corollary}
\label{JP.cor:cor3}
Let $0<A$, $1<p<\infty,$ $0<\eps <1/7\sqrt{2}$, $\alpha=49\eps^2/(1-49\eps^2)$,
and set
$d=\lfloor 4\bigl(A+(A/\eps)^{2/p}\bigr)^2\rfloor+1$.

If $\{e_k\}_{k=1}^{N} \subset L^2 (\R)$ is an orthonormal set such that 
for all $1 \leq k \leq N$,
$$
|\mu_p(e_k)|,\ |\mu_p(\widehat{e_k})|,\ 
\Delta_p(e_k),\ \Delta_p(\widehat{e_k})\leq A,
$$
then $$N \leq N_\C^s(\alpha,d)\leq N_\R^s(\alpha,2d).$$  
\end{corollary}

\begin{proof}
According to Lemma \ref{JP.lem:Po-prolate}, $e_1,\ldots,e_n$ 
are in $\pp_{A+(A/\eps)^{2/p},\eps}$.
The definition of $d$ and Theorem \ref{JP.th:LPth12} show that
$\{ e_j \}_{j=1}^N \subset \ss_{7\eps,d}$.
According to Proposition \ref{JP.prop:estimate},
$N\leq A_\C(7\eps,d)\leq N_\C^s(\alpha,d)\leq N_\R^s(\alpha,2d)$, where $\alpha=49\eps^2/(1-49\eps^2)$. 
\end{proof}


This approach does not, in general, give sharp results.  For example,
in the case $p=2$ the bound obtained by Corollary \ref{JP.cor:cor3} is not as good as the one
given in Section 2.  To see this, assume that $p=2$ and $A\geq 1$.
Then $4A^2(1+1/\eps)^2\leq d\leq 5A^2(1+1/\eps)^2$.
In order to apply the Delsarte, Goethals, Seidal bound (\ref{delsarte_bnd.ap}) 
we will now chose
$\eps$ so that $\alpha<\frac{1}{2\sqrt{d}}$ which will then give that $n\leq 4d$.
Our aim is thus to take $d$ is as small as possible by chosing
$\eps$ as large as possible.

For this, let us first take $\eps\leq 1/50$ so that $\alpha\leq 50\eps^2$.
It is then enough that $50\eps^2\leq\frac{1}{4A(1+1/\eps)}$, that is
$\eps^2+\eps-\frac{1}{200A}\leq 0$. We may thus take
$\eps= \frac{\sqrt{1+\frac{1}{50A}}-1}{2}$. Note that, as $A\geq 1$,
we get that $\eps\leq\frac{\sqrt{1+\frac{1}{50}}-1}{2}<1/50$.
This then gives
$$
n\leq 20d\leq 20A^2\left(1+\frac{2}{\sqrt{1+\frac{1}{50A}}-1}\right)^2
=20A^2\left(1+100A\Bigl(\sqrt{1+\frac{1}{50A}}+1\Bigr)\right)^2
\leq C A^4.
$$
In particular, the combinatorial methods allow one to take $N= CA^4$ in Corollary
\ref{JP.cor:cor3}, whereas the sharp methods of Section \ref{ap.sharp.md.sec} 
give $N = 8 \pi A^2,$ see Corollary \ref{ap.cor.opt.md}.

\subsection{The Quantitative Umbrella Theorem}
\label{JP.sec:QUT}
A second application of our method is a quantitative form of Shapiro's umbrella theorem.
As with the mean-dispersion theorem, Shapiro's proof does not provide a bound on the number 
of elements in the sequence.  As before, the combinatorial approach is well suited to this setting
whereas the approach of Section \ref{sec2.ap} is not easily applicable.

Given $f \in L^2 (\mathbb{R})$ and $\eps>0$, define
$$
C_f(\eps)=\inf\left\{ T \geq 0 \,:\ \int_{|t|>T}| f(t)|^2\leq\eps^2\norm{f}_2^2
\right\}.$$  Note that if $f$ is not identically zero then
for all $0<\eps<1$ one
has $0 < C_f (\eps) < \infty$.

\begin{theorem}
\label{JP.th:umbrella}
Let $\ffi,\psi\in L^2(\R)$ and 
$M = \min \{ \norm{\ffi}_2,\norm{\psi}_2 \} \geq 1$.
Fix $\frac{1}{50M} \geq \eps >0$, 
$T >\max \{ C_\ffi(\eps),C_\psi(\eps) \}$,
and $d = \lfloor 4T^2\rfloor+1$.  

If $\{e_n\}_{n=1}^N$ is an orthonormal sequence in $L^2(\R)$
such that for all $1 \leq n \leq N,$ and for almost all $x, \xi \in\R$, 
\begin{equation}
\label{JP.eq:umbrella}
|e_n(x)|\leq|\ffi(x)|\quad\mbox{and}\quad|\widehat{e}_n(\xi)|\leq|\psi(\xi)|,
\end{equation}
then 
\begin{equation} \label{ap.bnd.umbrella.thm}
N 
\leq N^s_\C(50\eps^2M^2,d) \leq N^s_\R(50\eps^2M^2,2d).
\end{equation}

In particular, $N$
is bounded by an absolute constant depending
only on $\ffi$ and $\psi$.
\end{theorem}

\begin{proof} 
By (\ref{JP.eq:umbrella}), $T>\max \{ C_\ffi(\eps),C_\psi(\eps) \}$,
implies $\{ e_n \}_{n=1}^N \subset\pp_{T,\eps M}$. 
According to Theorem \ref{JP.th:LPth12},
$\pp_{T,\eps M}\subset \ss_{7\eps M,d}$.
It now follows from Proposition \ref{JP.prop:estimate}, that
$$
N \leq \mathbb{A}_{\C} \left( 7 \eps M, d \right)
\leq N^s_\C \left( \frac{49\eps^2M^2}{1 - 49 \eps^2 M^2},d \right) 
\leq N^s_\C(50\eps^2M^2,d) \leq N^s_\R(50\eps^2M^2,2d)
.
$$
\end{proof}

Let us give two applications where one may get an explicit upper bound
by making a proper choice of $\eps$ in the proof above.

\begin{proposition}
\label{JP.ex:power}
Let $1/2<p$ and 
$\sqrt{\frac{2p-1}{2}} \leq C$ 
be fixed.
If $\{e_n \}_{n=1}^N \subset L^2(\R)$ is an orthonormal set such that
for all $1 \leq n \leq N$, and for almost every $x,\xi \in \R$,
$$|e_n(x)|\leq\frac{C}{(1+|x|)^p} \ \ \ \hbox{ and } \ \ \
|\widehat{e}_n(\xi)|\leq\frac{C}{(1+|\xi|)^p},$$
then 
$$
N\leq \begin{cases}
9^{\left(\frac{200\sqrt{2}C}{\sqrt{2p-1}}\right)^{\frac{4}{2p-1}}}
&\mbox{ if } 1/2 < p,\\
16\left(\frac{400C^2}{2p-1}\right)^{\frac{1}{p-1}}&\mbox{ if } 
1 < p \leq 3/2,\\
{4} 
\left(\frac{500C^2}{2p-1}\right)^{\frac{2}{2p-3}}&\mbox{ if } 3/2 < p.\\
\end{cases}$$

\end{proposition}

\begin{proof} 
If $\ffi(x)=\frac{C}{(1+|x|)^p}$, then 
$M = ||\ffi||_2 = C\sqrt{\frac{2}{2p-1}} \geq 1$, 
and a computation for $0 < \eps \leq 1$ shows that
$C_\ffi(\eps)=\frac{1}{\eps^{2/(2p-1)}}-1.$
Let $\delta =\delta(\eps)=\frac{4}{\eps^{4/(2p-1)}}$ and
$\alpha = \alpha(\eps) =\frac{100 C^2\eps^2}{2p-1}$.
Taking $T = C_{\ffi}(\eps)$ implies that 
$d = \lfloor 4T^2\rfloor+1
\leq \delta(\eps)$.

If $0 < \eps \leq \frac{1}{50M}$, then
Theorem \ref{JP.th:umbrella} gives the bound
$N\leq N^s_\C\left(\alpha (\eps ),\delta (\eps ) \right)$.
We shall chose $\eps$ differently for the various cases.

{\em Case 1.}  For the case $1/2 <p,$ take 
$\eps=\frac{1}{50M},$ 
and use the exponential bound given by
Lemma \ref{ap.expbnd.lem} for $N^s_\C(\alpha (\eps), \delta (\eps))$
to obtain the desired estimate.

{\em Case 2.}  For the case  $1<p\leq 3/2$, 
let $\eps_0=
\left(\frac{\sqrt{2p-1}}{20C}\right)^{\frac{2p-1}{2(p-1)}}$,
$\alpha=\alpha(\eps_0)$, and $\delta=\delta(\eps_0)$.
Note that $\alpha = \frac{1}{2\sqrt{\delta}}< \frac{1}{\sqrt{2\delta}}$, 
and also
that $\eps_0 \leq \frac{1}{50M},$ since $1 < p \leq 3/2$.
Thus the bound (\ref{delsarte_bnd.ap}) yields
$N \leq N^s_\R(\alpha,2\delta)= 4(1-\alpha^2)\delta\leq 4\delta$. 
The desired estimate follows.

{\em Case 3.}  For the case $3/2<p$, define 
$\eps_1=\left(\frac{\sqrt{2p-1}}{50C\sqrt{2}}\right)^{\frac{2p-1}{2p-3}}$
and note that  $\eps_1 \leq \frac{1}{50M}.$
Since $3/2<p$, taking
$\eps<\eps_1$, $\alpha=\alpha(\eps)$, $\delta=\delta(\eps),$ implies
that $\alpha(\eps) < 1/ \delta(\eps)$.
Thus, by Lemma \ref{ap.sphcodebnd.lem}, $N \leq \delta(\eps)$ 
for all $\eps < \eps_1$.  Hence,
$N \leq \delta (\eps_1),$ and the desired estimate follows.

\end{proof}

Note that in the case $1/2 <p$, the upper bound
in Proposition \ref{JP.ex:power} approaches infinity 
as $p$ approaches $1/2$.  Indeed, we refer the reader to
the counterexamples for $p<1/2$ in \cite{JP.DeR, By.ap}.
The case $p=1/2$ seems to be open as \cite{JP.DeR} need an extra logarithmic factor
in their construction. For perspective in the case $3/2<p$, 
if one takes $C = C_p = \sqrt{\frac{2p -1}{2}}$,
then the upper bound in Proposition \ref{JP.ex:power} 
approaches $4$ as $p$ approaches infinity.

\begin{proposition}
Let $0<a\leq 1$ and $(2a)^{1/4} \leq C$ be fixed.
If $\{e_n\}_{n=0}^N \subset L^2 (\R)$ is an orthonormal set such that
for all $n$ and for almost every $x, \xi \in \R$
$$|e_n(x)|\leq Ce^{-\pi a|x|^2} \ \ \ \hbox{ and } \ \ \
|\widehat{e}_n(\xi)|\leq Ce^{-\pi a|\xi|^2},$$
then 
$$N \leq  2+ \frac{8}{a\pi} \max \Big\{ 
2\ln \left( \frac{50 C \sqrt{\pi} e^{\pi}}{a^{1/4}} \right),
\ln \left( \frac{50 \pi C^2 e^{\pi a /2}}{a^{5/2}e^{2\pi}} \right)
\Big\}.$$
\end{proposition}

\begin{proof}  Let $\gamma_a(x)=Ce^{-\pi a|\xi|^2}$ and let
$C_a(\eps)=C_{\gamma_a}(\eps)$.
First note that
\begin{eqnarray*}
\int_{|t|>T} |\gamma_a (t)|^2 dt &=&
\int_{|t|>T}C^2e^{-2\pi a|t|^2} {d}t
=\frac{2C^2}{\sqrt{a}}\int_{T\sqrt{a}}^{\infty}
\frac{(1+s^2)e^{-2\pi s^2}}{1 +s^2} {d}s\\
&\leq&
\frac{C^2\pi (1+aT^2)}{\sqrt{a}} e^{-2\pi a T^2},
\end{eqnarray*}
while $M=||\gamma_a||_2 = 
(\int_{\R} C^2e^{-2\pi a|t|^2} {d}t)^{1/2}
=\frac{C}{(2a)^{1/4}}$. In particular, $\norm{\gamma_a}_2\geq1$.
Now for every $T>0$, set 
$\eps(T)=2^{-1/4} \sqrt{\pi} \sqrt{1+aT^2}e^{-\pi a T^2}$, so that
$C_a\bigl(\eps(T)\bigr) \leq T$.

By Theorem \ref{JP.th:umbrella}, we get that 
$N\leq N_\R^s(50\eps^2(T)M^2,8T^2+2)$, provided $\eps(T) \leq \frac{1}{50M}$.
Let us first see what condition should be imposed on $T$ to have 
$\eps(T) \leq \frac{1}{50M}$.  Setting $s=(1+aT^2)$, this condition is 
equivalent to $\sqrt{s} e^{- \pi s} \leq 
\frac{e^{-\pi} a^{1/4}}{50 C\sqrt{\pi}}.$
Thus, it suffices to take $s \geq  
\frac{2}{\pi} \ln \left( \frac{50 C \sqrt{\pi} e^{\pi}}{a^{1/4}} \right)$, 
and
$T^2 \geq \frac{2}{a\pi} 
\ln \left( \frac{50 C \sqrt{\pi} e^{\pi}}{a^{1/4}} \right)$.

We will now further choose $T$ large enough to have $50\eps(T)^2M^2<
\frac{1}{8T^2+2}$,
so that Lemma \ref{ap.sphcodebnd.lem} will imply 
$N \leq N_\R^s(50\eps(T)^2M^2,8T^2+2)=8T^2+2$.
This time, the condition reads
$(1 + aT^2)(1 +4T^2) e^{-2\pi a T^2} < \frac{\sqrt{a}}{50 \pi C^2}$.
Let $r= a(4T^2 +1)$.  Thus, it suffices to take
$r^2 e^{-\frac{\pi}{2} r} < 
\frac{a^{5/2} e^{2\pi}}{50 \pi C^2 e^{{\pi a}/{2}}}$
It is enough to take 
$r> \frac{4}{\pi} \ln \left( \frac{50 \pi C^2 
e^{\pi a /2}}{a^{5/2}e^{2\pi}} \right)$,
and $T^2 > \frac{1}{a\pi} 
\ln \left( \frac{50 \pi C^2 e^{\pi a /2}}{a^{5/2}e^{2\pi}} \right)$.

Combining the bounds for $T^2$ from the previous two paragraphs
yields
$$N \leq  2+ \frac{8}{a\pi} \max \Big\{ 
2\ln \left( \frac{50 C \sqrt{\pi} e^{\pi}}{a^{1/4}} \right),
\ln \left( \frac{50 \pi C^2 e^{\pi a /2}}{a^{5/2}e^{2\pi}} \right)
\Big\}.$$


\end{proof}


A careful reading of the proof of the Umbrella Theorem shows the following:
\begin{proposition}
Let $0<C,$ and let $1\leq p,q,\widehat{p},\widehat{q}\leq \infty$ 
satisfy $\dst\frac{1}{p}+\frac{1}{q}=1$ and 
$\dst\frac{1}{\widehat{p}}+\frac{1}{\widehat{q}}=1$. Let $\ffi\in L^{2p}(\R)$ and $\psi\in L^{2\widehat{p}}(\R)$, and suppose that 
$\ffi_k\in L^{2q}(\R)$ and $\psi_k\in L^{2\widehat{q}}(\R)$ 
satisfy $\norm{\ffi_k}_{2q} \leq C$, $\norm{\psi_k}_{2\widehat{q}}\leq C$.
There exists $N$ such that, if $\{e_k\}\subset L^2 (\R)$ 
is an orthonormal set which for all $k$ and
almost every $x, \xi \in\R$ satisfies
$$
|e_k(x)|\leq \ffi_k(x) \ \ffi(x) \quad \ \ \ \hbox{ and } \ \ \
\quad |\widehat{e_k} (\xi)|\leq \psi_k(\xi) \ \psi(\xi),
$$
then $\{e_k\}$ has at most $N$ elements.  As with previous results, 
a bound for $N$ can be obtained in terms of spherical codes.  The
bound for $N$ depends only on $\ffi, \psi, C$.
\end{proposition}

Indeed, let $\eps>0$ and take $T>0$ big enough to have $\dst\int_{|t|>T}|f(t)|^{2p}\mbox{d}t\leq\eps^p/C^{2/p}$. Then
\begin{eqnarray*}
\int_{|t|>T}|e_k(t)|^2\,\mbox{d}t
&\leq&\int_{|t|>T}|\ffi_k(t)\ffi(t)|^2\,\mbox{d}t
\leq \left(\int_{|t|>T}|\ffi_k(t)|^{2q}\,\mbox{d}t\right)^{1/q}\left(\int_{|t|>T}|\ffi(t)|^{2p}\,\mbox{d}t\right)^{1/p}\\
&\leq& C^2(\eps^p/C^{2/p})^{1/p}=\eps.
\end{eqnarray*}
A similar estimate holds for $\widehat{e_k}$ and we conclude as in the proof of the Umbrella Theorem.

\subsection{Angles in Riesz bases} \label{ap.rieszang.sec}

Let us now conclude this section with a few remarks on Riesz bases.
Let $\{x_k\}_{k=0}^{\infty}$  be a Riesz basis for $L^2(\R)$ with orthogonalizer $U$
and recall that, for every sequence $\{a_n\}_{n=0}^{\infty} \in\ell^2$,
\begin{equation}
\frac{1}{\norm{U}^2}\sum_{n=0}^{\infty}|a_n|^2
\leq \norm{\sum_{n= 0}^{\infty} a_nx_n}_2^2\leq\norm{U^{-1}}^2\sum_{n=0}^{\infty} |a_n|^2.
\label{JP.eq:riesz2}
\end{equation}

Taking $a_n=\delta_{n,k}$ in (\ref{JP.eq:riesz2}) shows that 
$\frac{1}{\norm{U}}\leq\norm{x_k}_2 \leq\norm{U^{-1}}$.
Then taking $a_n=\delta_{n,k}+\lambda\delta_{n,l}$, $k\not=l$ and $\lambda=t,-t$, $t>0$ gives
$$
\frac{1}{\norm{U}^2}(1+t^2)\leq\norm{x_k}_2^2+t^2\norm{x_l}_2^2+2t\abs{\mbox{Re}\scal{x_k,x_l}}
\leq\norm{U^{-1}}^2(1+t^2)
$$
thus $\abs{\mbox{Re}\scal{x_k,x_l}}^2$ is 
\begin{eqnarray*}
&\leq&\min\Bigl((\norm{x_k}_2^2-\norm{U}^{-2})(\norm{x_l}_2^2-\norm{U}^{-2}),
(\norm{U^{-1}}^2-\norm{x_k}_2^2)(\norm{U^{-1}}^2-\norm{x_l}_2^2)\Bigr)\\
&\leq&\norm{x_k}_2^2\norm{x_l}_2^2\min\ent{
\left(1-\frac{1}{\norm{U}^2\norm{U^{-1}}^2}\right)^2
,\left(\frac{\norm{U^{-1}}^2}{\norm{U}^2}-1\right)^2}
\end{eqnarray*}
while taking $\lambda=it,-it$, $t>0$ gives the same bound for $\abs{\mbox{Im}\scal{x_k,x_l}}^2$.
It follows that
\begin{equation}
\label{JP.eq:c(u)}
\abs{\scal{x_k,x_l}}\leq C(U)\norm{x_k}_2 \norm{x_l}_2 \leq C(U)\norm{U^{-1}}^2
\end{equation}
where
$$
C(U):=\sqrt{2}
\min\ent{1-\left(\frac{1}{\norm{U}\norm{U^{-1}}}\right)^2,\left(\frac{\norm{U^{-1}}}{\norm{U}}\right)^2-1}.
$$
We may now adapt the proof of Proposition \ref{JP.prop:estimate} to Riesz basis:

\begin{proposition}
\label{JP.prop:estimateriesz}
Let $\{\psi_k)_{k=1}^{\infty}$ be an orthonormal basis for $L^2(\R)$.
Fix $d\geq 0$ and let $\P_d$ be the projection on the span of $\{\psi_1,\ldots,\psi_{d-1}\}$.

Let $\{x_k\}_{k=1}^{\infty}$ be a Riesz basis for $L^2(\R)$ and let $U$ be its orthogonalizer. Let $\eps>0$ be such that
$\eps<\min\left(\norm{U}^{-2},\sqrt{\frac{\norm{U}^{-2}-C(U)\norm{U^{-1}}^2}{2}}\right)$ and let
\begin{equation}
\label{JP.eq:defalpha}
\alpha=\frac{\eps^2+C(U)\norm{U^{-1}}^2}{\norm{U}^{-2}-\eps^2}.
\end{equation}
If $\{x_k\}_{k=1}^N$ satisfies
$\norm{x_k-\P_dx_k}_2<\eps$ then $N \leq N_\K^s(\alpha,d)$.
\end{proposition}

\begin{proof}
Assume without loss of generality that $x_0,\ldots,x_N$ satisfy $\norm{x_k-\P_dx_k}<\eps$
and let $a_{k,j}=\scal{x_k,\psi_j}$.

Write $v_k=(a_{k,1},\ldots,a_{k,d})\in\K^d$ then, the same computation as in (\ref{JP.eq:orsp}), for $k\not=l$
$$
\scal{v_k,v_l}=\scal{x_k-\P_dx_k,\P_dx_l-x_l}+\scal{x_k,x_l}
$$
thus $\abs{\scal{v_k,v_l}}\leq\eps^2+\abs{\scal{x_k,x_l}}$.
On the other hand
$$
\norm{v_k}=\norm{\P_dx_k}=(\norm{x_k}^2-\norm{x_k-\P_df_k}^2)^{1/2}\geq(\norm{U}^{-2}-\eps^2)^{1/2}
$$
It follows from (\ref{JP.eq:c(u)}) that $w_k=\dst\frac{v_k}{\norm{v_k}}$ satisfies, for $k\not=l$,
$$
\abs{\scal{w_k,w_l}}\leq\frac{\eps^2+C(U)\norm{U^{-1}}^2}{\norm{U}^{-2}-\eps^2}
$$
and $\{w_k\}$ forms a spherical $[-\alpha,\alpha]$-code in $\K^d$.
\end{proof}

Note that the condition on $\eps$ implies that $0<\alpha<1$.
Also note that if $U$ is a near isometry in the sense that $(1+\beta)^{-1}\leq\norm{U}^2\leq\norm{U^{-1}}^2\leq 1+\beta$
then $C(U)\leq\sqrt{2}\frac{\beta(2+\beta)}{(1+\beta)^2}$ and
$\alpha\leq\frac{(1+\beta)\eps^2+\beta(2+\beta)}{1-(1+\beta)\eps^2}$. In particular, if $U$ is near enough to an isometry,
meaning that $\beta$ is small enough, then this $\alpha$ is comparable with the $\alpha$ of Proposition
\ref{JP.prop:estimate}.

As a consequence, we may then easily adapt the proof of results that relied on Proposition \ref{JP.prop:estimate} 
to the statements about Riesz bases. For example, an Umbrella Theorem for Riesz bases
reads as follows:

\begin{theorem}
\label{JP.th:umbrellaR}
Let $\ffi,\psi\in L^2(\R)$ with $\norm{\ffi}_2,\norm{\psi}_2 \geq1$.
Let $\{f_n\}_{n=1}^{\infty}$ 
be a Riesz basis for $L^2(\R)$ 
with orthonormalizer $U$ that is near enough to an isometry
$(1+\beta)^{-1}\leq\norm{U}^2\leq\norm{U^{-1}}^2\leq 1+\beta$ with $\beta$ small enough.
Then there exists a constant $N=N(\ffi,\psi,\beta)$ depending only on $\ffi, \psi$ and 
$\beta$, such that the number of terms of the basis that satisfies
$$|f_n(x)|\leq|\ffi(x)|\quad\mbox{and}\quad|\widehat{f}_n(\xi)|\leq|\psi(\xi)|$$
for almost all $x, \xi\in\R$ is bounded by $N$.  As with previous results,
a bound on $N$ can be given in terms of spherical codes.
\end{theorem}

\subsection*{Acknowledgements}  A portion of this work was performed during the Erwin Schr\"odinger
Institute (ESI) Special Semester on ``Modern methods of time-frequency analysis.''
The authors gratefully acknowledge ESI for its hospitality
and financial support.  The authors also thank 
Professor H.S. Shapiro for valuable comments related to the material.

\end{document}